\documentclass[10pt]{article}
\usepackage{latexsym}\usepackage{amsbsy}\usepackage{amssymb}
\usepackage[latin1]{inputenc} 

\newcommand{\D}{\ensuremath{{\mathcal D}} } 
  \newcommand{\beqs}{\begin{displaymath}}
\newcommand{\eeqs}{\end{displaymath}}

\newcommand{\ep}{\hspace*{\fill}$\Box$} \newcommand{\eps}{\varepsilon}
 
 \newcommand{\R}{\mathbb R}
\newcommand{\N}{\mathbb N} \newcommand{\C}{\mathbb C}
 \newcommand{\K}{\mathbb K}
\newcommand{\gK}{{\cal K}} \newcommand{\gR}{{\cal R}}
\newcommand{\gC}{{\cal C}} 
\newcommand{\gs}{\ensuremath{{\mathcal G}} }

\newcommand{\es}{\ensuremath{{\mathcal E}} }
\newcommand{\esm}{\ensuremath{{\mathcal E}_M} }

\newcommand{\ns}{\ensuremath{{\mathcal N}} }
\newcommand{\ks}{\ensuremath{{\mathcal K}} }

\newcommand{\comp}{\subset\subset} \newcommand{\cinfty}{{\cal C}^\infty}
\newtheorem{thr}{\hspace*{-3mm} \bf}[section] \newcommand{\bt}{\begin{thr}
{\bf Theorem. }} \newcommand{\et}{\end{thr}} \newcommand{\bp}{\begin{thr}
{\bf Proposition. }} \newcommand{\bc}{\begin{thr} {\bf Corollary. }}
\newcommand{\blem}{\begin{thr} {\bf Lemma. }} \newcommand{\bex}{\begin{thr}
{\bf Example. }\rm} \newcommand{\bexs}{\begin{thr} {\bf Examples. }\rm}
\newcommand{\bd}{\begin{thr} {\bf Definition. }}
\newcommand{\beast}{\begin{eqnarray*}} \newcommand{\eeast}{\end{eqnarray*}}
 \newcommand{\al}{\alpha}
 \newcommand{\ga}{\gamma}
\newcommand{\Om}{\Omega}\newcommand{\Ga}{\Gamma}
 
\newcommand{\vphi}{\varphi}

\newcommand{\Vol}{\mbox{Vol\,}}

 \newcommand{\G}{{\cal G}}
 \newcommand{\CC}{{\cal C}}    

 \newcommand{\beq}{ \begin{equation} }
\newcommand{\eeq}{\end{equation} }
\newcommand{\bea}{\begin{eqnarray}}\newcommand{\eea}{\end{eqnarray}}
\newcommand{\beas}{\begin{eqnarray*}}\newcommand{\eeas}{\end{eqnarray*}}

\newcommand{\ben}{\begin{enumerate}}\newcommand{\een}{\end{enumerate}}
\newcommand{\ba}{\begin{array}}\newcommand{\ea}{\end{array}}
\newcommand{\brem}{\begin{thr} {\bf Remark. }\rm}
\newcommand{\ethi}{\end{thr}} 

\newcommand{\gsh}{\ensuremath{{\mathcal G}^{\mathit h}} }

\newcommand{\cc}{{\mathcal C}}

\begin{document}

\title{Non-smooth differential geometry and algebras of generalized functions}

\author{Michael Kunzinger \footnote{Electronic mail:
        michael.kunzinger@univie.ac.at} }

        \date{\small Dedicated to John Horv\'ath on the occasion of his 80th birthday. \\ 
                  With special thanks for all his support over the years.} \maketitle

\begin{abstract}
Algebras of generalized functions offer possibilities beyond the purely distributional
approach in modelling singular quantities in non-smooth differential geometry. This
article presents an introductory survey of recent developments in this field and
highlights some applications in mathematical physics.

 
\end{abstract}

\section{Introduction}\label{intro}
Non-smooth differential geometry provides an important tool in a variety of applications,
in particular in mathematical physics. As examples we mention non-smooth Hamiltonian
mechanics (\cite{marsden,marsden2}) and the analysis of singular spacetimes in
general relativity (cf., e.g., \cite{gt}, \cite{herbertgeo}, \cite{wilson}, and
\cite{vickersESI} for a recent survey). Linear distributional geometry (\cite{deR,
marsden, parker}) is only of limited use in a genuinely nonlinear context, as,
e.g., in general relativity, where the nonlinearity of the Einstein field
equations and the interest in curvature quantities introduces requirements on
the underlying theory of generalized functions which distribution theory
is unable to meet. A nonlinear extension of linear distributional geometry 
displaying promising capabilities for overcoming these conceptual problems
has been developed over the past years based on Colombeau's theory of generalized
functions. It is the aim of the present paper to provide an introduction to
this field and some of its applications.

In the remainder of this section we fix some notation and terminology from
differential geometry and distribution theory. Section \ref{col} gives a
quick introduction to some of the fundamental ideas of Colombeau theory 
both in the local and in the manifold setting. In Section \ref{mfval} 
we consider generalized functions taking values in a
differentiable manifold, a construction which has no analogue in distribution
theory yet is of central importance for nonlinear distributional geometry
as it allows to formulate a functorial theory of generalized functions
in a global context. In particular, it allows to introduce notions like
flows of generalized vector fields or geodesics of generalized metrics.
Finally, in Section \ref{riemann} we develop a generalized pseudo-Riemannian
geometry in this setting and give some applications of the resulting theory
in general relativity.

In what follows, $X$ and $Y$ always mean paracompact, smooth Hausdorff manifolds 
of dimension $n$ resp.\ $m$. We denote vector bundles with base space $X$ by
$(E,X,\pi_X)$ or $E\to X$ for short and write a vector bundle chart over 
the chart $(V,\psi)$ of $X$ as $(V,\Psi)$.
For vector bundles $E\to X$ and $F\to Y$, by $\mathrm{Hom}(E,F)$
we mean the space of vector bundle homomorphisms from $E$ to $F$.
Given $f\in\mathrm{Hom}(E,F)$ the unique smooth map from $X$ to $Y$
satisfying $\pi_Y\circ f = \underline{f}\circ \pi_X$
is denoted by $\underline{f}$.  
For vector bundle charts $(V,\Phi)$ of $E$ and $(W,\Psi)$ of $F$ we
write the local vector bundle homomorphism
$f_{\Psi\Phi}:=
\Psi\circ f \circ \Phi^{-1}: \vphi(V\cap \underline{f}^{-1}(W))
\times \K^{n'} \to \phi(W) \times \K^{m'}$
in the form
\[
f_{\Psi\Phi}(x,\xi) =
(f_{\Psi\Phi}^{(1)}(x),f_{\Psi\Phi}^{(2)}(x)\cdot\xi)\,.
\]

The space of smooth sections of a vector bundle $E\to X$
is denoted by $\Ga(X,E)$. $T^r_s(X)$ is the $(r,s)$-tensor bundle over $X$ 
and we use the following notation for spaces of tensor
fields  ${\mathcal T}^r_s(X):=\Ga(X,T^r_s(X))$, ${\mathfrak X}(X):=\Ga(X,TX)$ and
${\Om}^1(X):=\Ga(X,T^*X)$, where $TX$ and $T^*X$ denote the tangent
and cotangent bundle of $X$, respectively. ${\cal P}(X,E)$ is
the space of linear differential operators $\Gamma(X,E) \to \Gamma(X,E)$. 
For $E=X\times \R$ we write ${\cal P}(X)$ instead of ${\cal P}(X,E)$. 

We denote by $\Vol(X)$ the volume bundle over $X$, its smooth sections 
are called one-densities. The space $\D'(X,E)$ of $E$-valued distributions 
on $X$ is defined as the topological 
dual of the space of compactly supported sections of the bundle $E^*\otimes\Vol(X)$:
$$
  \D'\,(X,E):=[\Ga_c(X,E^*\otimes\Vol(X))]'\,.
$$ For $E=X\times \R$ we obtain $\D'(X):=\D'(X,E)$, the space of distributions on $X$.
The isomorphism of $\CC^\infty(X)$-modules
$$
\D'(X,E) \, \cong \,
   \D'(X)\otimes_{\CC^\infty(X)}\Ga(X,E)
$$
shows that $E$-valued distributions can be viewed as sections with distributional 
coefficients.

\section{Colombeau generalized functions on differentiable manifolds} \label{col}
When trying to extend linear distribution theory to a nonlinear theory
of generalized functions one is faced with certain fundamental obstacles.
To give a simple example, let $vp \frac{1}{x}$ be the Cauchy principal
value of $1/x$ on $\R$. Then since 
$$
0 = (\delta(x)\cdot x)\cdot vp\frac{1}{x} \not= \delta(x)\cdot(x\cdot 
vp\frac{1}{x}) = \delta(x)
$$
it follows that the usual multiplication on $C^\infty\times \D'$ cannot
be extended to an associative and commutative multiplication on $\D'\times
\D'$. Similarly, it can be shown that $\D'$ cannot be endowed with the 
structure of an associative commutative algebra compatible with the usual
product in $L^\infty$: with $H$ the Heaviside function, the fact that
$H^2 = H$ would by the Leibniz rule entail $(H^2)'=2HH'$, $(H^3)'=3H^2H'$,
so $2HH'=H'=3HH'$. But then $\delta = H' = 0$, a contradiction. For
a comprehensive analysis of the problem of multiplication of distributions
see \cite{MObook}. 

Apart from nonlinear analysis on certain (function-) subalgebras of $\D'$
(Sobo\-lev spaces) the second main option therefore consists in embedding the space 
of distributions into an appropriate (associative and commutative) algebra $\gs$ 
of generalized functions, the aim being to retain as many of the standard features 
of distribution theory as possible. In particular, we want $\gs$ to be a
differential algebra with unit $f(x)\equiv 1$ and derivation operators
extending those on $\D'$. Our previous example demonstrates
that under these assumptions the product in $\gs$ cannot extend the pointwise
product of functions in $L_{loc}^\infty$. Furthermore, by a celebrated result
of L.\ Schwartz (\cite{Schw}), it cannot extend the pointwise product 
of $\cc^k$-functions for any $k\in \N_0$ either. Due to these differential-algebraic
constraints the maximal possible compatibility of the product $\cdot$ in $\gs$
is that $\cdot\mid_{\cc^\infty\times \cc^\infty}$ coincide with the usual pointwise
product of functions.

Differential algebras satisfying this maximal set of requirements were first
constructed by J.F.\ Colombeau in the early 1980ies (\cite{CPort,CJMAA,CCR,c1,c2}). 
The basic principles underlying his approach are regularization through convolution
and asymptotic estimates in terms of a regularization parameter. In the so-called
special version of the construction, $\D'(\R^n)$ is embedded into a certain subalgebra 
$\esm(\R^n)$ of $\cc^\infty(\R^n)^I$ (with $I:=(0,1]$) through convolution:
$$
\D'(\R^n) \ni w \mapsto (w*\rho_\eps)_{\eps\in I}\,.
$$
Here $\rho$ is a Schwartz function with $\int \rho = 1$ and $\rho_\eps(x)=1/\eps^n\rho(x/\eps)$.
$\cc^\infty(\R^n)^I$ is a differential algebra with operations defined componentwise and
the above map is obviously linear and commutes with partial derivatives. On the
other hand, a natural way of embedding $\cc^\infty(\R^n)$ into $\cc^\infty(\R^n)^I$ is
the diagonal embedding
$$
\cc^\infty(\R^n) \ni f \mapsto (f)_{\eps\in I}\,.
$$
Clearly this map preserves the pointwise product of smooth functions. The idea, 
therefore, is to factor $\esm(\R^n)$ by an ideal $\ns(\R^n)$ containing $(f*\rho_\eps 
- f)_\eps$ for each $f\in \cc^\infty(\R^n)$. The resulting quotient algebra would then 
satisfy the above maximal set of requirements on a differential algebra containing the
space of distributions. Now (assuming $n=1$ for the moment), Taylor's theorem gives
\begin{eqnarray*}
&&(f\ast \rho_\eps-f)(x)=\int\!(f(x-y)-f(x))\rho_\eps(y)\, dy\\
&&=\int\!\sum_{k=1}^m \frac{(-\eps y)^k}{k!}f^{(k)}(x)\rho(y)\, dy +
  \!\int\!\frac{(-\eps y)^{m+1}}{(m+1)!}f^{(m+1)}(x-\theta\eps y)\rho(y)\, dy.
\end{eqnarray*}
If we additionally suppose that $\int \rho(x)x^k\,dx=0$ for all $k\ge 1$ then this
expression converges to zero, faster than any power of $\eps$, uniformly on each compact 
set, in each derivative. The natural candidate for $\ns(\R^n)$ therefore is
\begin{eqnarray*}
&&\ns(\R^n) = \{(u_\eps)_\eps\in \cc^\infty(\R^n)^I \mid \forall K \comp \R^n\  
\forall \alpha\in \N_0^n\ \forall m \in \N \mbox{ : }\\
&& \hphantom{\ns(\R^n) =\{(u_\eps)_\eps\in \cc^\infty(\R^n)^I \mid}
\sup_{x\in K}|\partial^\alpha u_\eps(x)| = O(\eps^{m}) \mbox{ as } \eps\to 0 \}
\end{eqnarray*}
Elements of $\ns(\R^n)$ are called {\em negligible}. The definition of $\ns(\R^n)$
in turn fixes the maximal subalgebra $\esm(X)$ (the algebra of {\em moderate} nets) 
of $\cc^\infty(\R^n)^I$ in which $\ns(\R^n)$ is an ideal as
\begin{eqnarray*}
&&\esm(\R^n)=\{(u_\eps)_\eps\in \cc^\infty(\R^n)^I \mid \forall K \comp \R^n\  
\forall \alpha\in \N_0^n\ \exists N \in \N \mbox{ with }\\
&& \hphantom{\esm(\R^n)=\{(u_\eps)_\eps\in \cc^\infty(\R^n)^I \mid}
\sup_{x\in K}|\partial^\alpha u_\eps(x)| = O(\eps^{-N}) \mbox{ as } \eps\to 0 \}
\end{eqnarray*}
The (special) Colombeau algebra on $\R^n$ is then defined as the factor algebra
$\gs(\R^n)=\esm(\R^n)/\ns(\R^n)$. As indicated above, the map $\iota: \D'(\R^n) \to
\gs(\R^n)$, $\iota(w) =$ $[$class of $(w*\rho_\eps)_\eps]$ provides a linear embedding 
which coincides with the diagonal embedding
$\sigma: \cc^\infty(\R^n)\to \gs(\R^n)$, $\sigma(f)=$ $[$class of $(f)_\eps]$ on 
$\cc^\infty(\R^n)$, hence verifies all the requirements made above. From here one
may proceed, using partitions of unity and suitable cut-off functions to construct
embeddings $\D'(\Om)\hookrightarrow \gs(\Om)$ for any open subset $\Om$ of $\R^n$.
Instead, we turn directly to the manifold case (\cite{AB,RD,hermannbook,ndg}). 
The basic features of the following
definition are in close correspondence to the Euclidean case discussed above.

\bd Let $X$ be a smooth, paracompact Hausdorff manifold and set $\es(X):=(\cc^\infty(X))^ I$.
The Colombeau algebra $\gs(X)$ on $X$ is defined as the quotient $\esm(X)/\ns(X)$, where
\begin{eqnarray*}
&&\esm(X)=\{(u_\eps)_\eps\in \es(X) \mid \forall K \comp X\  
\forall P\in {\mathcal P}(X)\ \exists N \in \N \mbox{ : }\\
&& \hphantom{\esm(X)=\{(u_\eps)_\eps\in \es(X) \mid \qquad\qquad}
\sup_{p\in K}|Pu_\eps(p)| = O(\eps^{-N}) \mbox{ as } \eps\to 0 \}\\
&&\ns(X)=\{(u_\eps)_\eps\in \es(X) \mid \forall K \comp X\  
\forall P\in {\mathcal P}(X)\ \forall m \in \N \mbox{ : }\\
&& \hphantom{\ns(X)=\{(u_\eps)_\eps\in \es(X) \mid \ \qquad\qquad}
\sup_{p\in K}|Pu_\eps(p)| = O(\eps^m) \mbox{ as } \eps\to 0 \}
\end{eqnarray*}
\et
We write $u=[(u_\eps)_\eps]$ for the class of $(u_\eps)_\eps$ in $\gs(X)$.
Restrictions of elements of $\gs(X)$ to open subsets of $X$ are defined
componentwise on representatives and $\gs(\_)$ is seen to be a fine and
supple (but not flabby) sheaf of differential algebras (\cite{RD,DP,OPS}).

Our first fundamental observation concerning the structure of $\gs(X)$ is that 
$\ns(X)$ can be characterized as a subspace of $\esm(X)$ without resorting 
to derivatives (\cite{found}, Th.\ 13.1, \cite{ndg}, Sec.\ 4):
\begin{equation}\label{mgth}
\ns(X)=\{(u_\eps)_\eps\in \esm(X) \mid \forall K \comp X \ \forall m\in \N 
\mbox{ : } \sup_{p\in K}|u_\eps(p)| = O(\eps^m)\}
\end{equation}
This characterization is a very convenient means both within Colombeau theory 
(as we shall see shortly) and in applications to partial differential equations 
(where it considerably simplifies uniqueness proofs). 

An important feature distinguishing Colombeau algebras from spaces of distributions
is the availability of a point value description of Colombeau functions. 
Componentwise insertion of points of $X$ into elements of $\gs(X)$ yields well-defined
{\em generalized numbers}, i.e., elements of the ring of constants $\gK:=\esm/\ns$
(with $\gK=\gR$ or $\gK=\gC$ for $\K=\R$ or $\K=\C$), where 
\begin{eqnarray*}
&& \esm = \{(r_\eps)_\eps \in \K^I \mid \exists N\in \N\ \mbox{: }|r_\eps| = O(\eps^{-N})\}\\
&& \ns \ =\ \! \{(r_\eps)_\eps \in \K^I \mid \forall m\in \N\ \mbox{: } |r_\eps| = O(\eps^{m})\}
\end{eqnarray*}
\bex \label{pvex}
Let $\vphi\in \D(\R)$, $\int \vphi = 1$, $\vphi_\eps(x):=\eps^{-1}\vphi(x/\eps)$ and
set $u_\eps(x):=\vphi_\eps(x-\eps)$. Then $u_\eps\to \delta$ in $\D'(\R)$, so $u:=[(u_\eps)_\eps]$
is not $0$ in $\gs(\R)$. Nevertheless, it is easily seen that every point value of every derivative 
of $u$ is zero in $\gK$.
\et
Thus point values on ``classical'' points $p\in X$ do not characterize elements of $\gs(X)$.
As can be seen in the above example, the reason for this failure is that Colombeau functions
are capable of modelling infinitesimal quantities which standard points are unable to detect.
Borrowing an idea from nonstandard analysis, the plan is therefore to introduce ``nonstandard
points'' which themselves may move around in the manifold in order to keep track of the
infinitesimal behavior of elements of $\gs(X)$. To this end we define an equivalence relation
$\sim$ on the space $X_c := \{(p_\eps)_\eps\in X^I\mid \exists K\comp X \ \exists \eps_0>0 \mbox{ s.t. }
p_\eps\in K \ \forall \eps<\eps_0\}$ as follows: for any Riemannian metric $h$ on $X$ with
distance function $d_h$, two nets $(p_\eps)_\eps$, $(q_\eps)_\eps$ are called equivalent,
$(p_\eps)_\eps\sim (q_\eps)_\eps$ if $d_h(p_\eps,q_\eps)=O(\eps^m)$ for each $m\in \N$.
We call $\tilde X_c:=X_c/\sim$ the space of compactly supported generalized points. Obviously
this definition does not depend on the specific Riemannian metric $h$. Then we have
\bt Let $u\in \gs(X)$ and $\tilde p = [(p_\eps)_\eps]\in \tilde X_c$. Then $u(\tilde p):=
[(u_\eps(p_\eps))_\eps]$ is a well-defined element of $\gK$. Moreover, $u=0$ if and only if
$u(\tilde p)=0$ in $\gK$ for all $\tilde p$ in $\tilde X_c$.
\et 
For the proof, see \cite{point, ndg}. To give an idea of the argument, let us have a look
at the case $X=\R^n$ (following \cite{OPS}, Prop.\ 3.1). 
If $u= 0 \in \gs(\R^n)$ and $p_\eps\in K\comp \R^n$ for $\eps$ small
then it is immediate from the definition of $\ns(\R^n)$ that $(u_\eps(p_\eps))_\eps\in \ns$,
i.e., $u(\tilde p) = 0 \in \gK$. Conversely, suppose that $u(\tilde p)=0$ for all $\tilde p\in
\tilde\R^n_c$ and let $K\comp \R^n$. For each $\eps\in I$ denote by $p_\eps$ the point in
$K$ where $|u_\eps|$ attains its maximum. Since $\tilde p = [(p_\eps)_\eps]\in \tilde\R^n_c$,
the negligibility estimates of order $0$ for $(u_\eps)_\eps$ on $K$ follow from
$(u_\eps(p_\eps))_\eps \in \ns$. But then
$u=0$ due to (\ref{mgth}).\ep\medskip

Note that in Example \ref{pvex}, $u(\tilde p)\not=0$ for $\tilde p=[(\eps)_\eps]$ if 
$\vphi(0)\not=0$. 

There are essentially two ways of connecting linear distribution spaces with Colombeau
algebras. Firstly, one can construct injective sheaf morphisms $\iota: \D'(\_)\hookrightarrow
\gs(\_)$. This can be done either using de Rham regularizations or, which basically amounts to the
same, directly by convolution with a fixed mollifier in charts (cf.\ \cite{RD, ndg}).
The resulting embedding is non-canonical, i.e.\ it depends on the ingredients of the 
construction (partition of unity, mollifier, cut-off functions, etc.). 
The main field of application of the special version of Colombeau
algebras therefore lies in areas where a regularization procedure for the singular
quantities to be modelled suggests itself by the nature of the problem (cf.\ \cite{MObook,
RD,book}).
For so-called full variants of Colombeau algebras on manifolds, allowing for a {\em canonical}
embedding of the space of distributions we refer to \cite{found,vim}.

The second link to linear distribution theory is the concept of association: two elements
$u$, $v$ of $\gs(X)$ are called {\em associated}, $u\approx v$ if $u_\eps - v_\eps \to 0$ in $\D'(X)$.
If $\int u_\eps \mu \to \langle w,\mu\rangle$ for some $w\in\D'(X)$ and each compactly supported
one density $\mu$, i.e., if $u_\eps \to w$ in $\D'(X)$ then $w$ is called associated
distribution to $u$. Clearly these definitions do not depend on the chosen representatives. Besides
this concept of ``equality in the sense of distributions'' one may also introduce more
restrictive equivalence relations on $\gs(X)$. In particular, we mention the concept
of $\cc^k$-association: $u,\ v\in \gs(X)$ are called $\cc^k$-associated, $u\approx_k v$ if for all
$l\le k$ and all $\xi_1,\dots,\xi_l \in {\frak X}(X)$, $L_{\xi_1}\dots L_{\xi_l} (u_\eps-v_\eps)
\to 0$, uniformly on compact sets. In applications it is often the case that modelling of
singular quantities and analytical treatment of the problem at hand (e.g., solution of a
nonlinear PDE) is carried out in $\gs$, while a distributional interpretation of the result 
is effected through the notion of association. Concerning the examples inspected at the
beginning of this section we note that, in $\gs(\R)$, $x\cdot \delta$
is associated but not equal to $0$  and $H^m \not= H$, but $H^m\approx H$ for all $m\in \N$.
This complies with the intuitive feeling that over and above the distributional picture,
modelling in $\gs$ allows to fix the ``microstructure'' of singular quantities, reflected
in a notion of equality which is more restrictive than equality in the distributional sense.
It can also be viewed as a further nonstandard aspect of the theory (cf.\ \cite{MObook}, §10 
for an in-depth discussion). 

For a vector bundle $E\to X$ we define the spaces of moderate resp.\ negligible sections as
\beas
        \Ga_{\esm}(X,E)&=& \{ (s_\eps)_{\eps\in I}\in \Ga(X,E)^I : 
                \ \forall P\in {\cal P}(X,E)\,\\
                 && \hspace{1cm}
                 \forall K\comp X \, \exists N\in \N:\ \sup_{p\in K}\|Pu_\eps(p)\| = O(\eps^{-N})\}\\ [.5em]
        \Ga_{\ns}(X,E)&=& \{ (s_\eps)_{\eps\in I}\in \Ga(X,E)^I : 
                \ \forall P\in {\cal P}(X,E)\, \\
                 &&\hspace{1.3cm}
                 \forall K\comp X \, \forall m\in \N:\ \sup_{p\in K}\|Pu_\eps(p)\| = O(\eps^{m})\}\\ 
\eeas
where $\|\ \|$ denotes the norm induced on the fibers of $E$ by any Riemannian metric. $\Ga_{\esm}(X,E)$
is a $\gs(X)$-module with submodule $\Ga_{\ns}(X,E)$ and we define the $\gs(X)$-module $\Ga_\gs(X,E)$ of generalized
sections of the bundle $E\to X$ as the quotient $\Ga_{\esm}(X,E)/\Ga_{\ns}(X,E)$. As in the scalar case we may omit all 
differential operators from the definition of $\Ga_{\ns}(X,E)$ if we suppose the $(s_\eps)_\eps$ 
to be moderate. Important special cases are the space $\gs^r_s(X)$ of generalized $(r,s)$-tensor fields
and the space $\bigwedge^k_{\gs}(X)$ of generalized $k$-forms, corresponding to 
$E=T^r_s(X)$ and $E=\bigwedge^k T^*X$, respectively. 

$\Ga_\gs(\_\, ,E)$ is a fine sheaf of $\gs(\_)$-modules. Its algebraic structure is 
clarified by the following theorem (\cite{ndg}, Sec.\ 6):
\bt The $\gs(X)$-module $\Ga_\gs(X,E)$ is projective and finitely generated. Moreover,
the following isomorphisms of $\cc^\infty(X)$-modules hold:
$$
\Ga_{\gs}(X,E) \cong \G(X)\otimes_{{\cal C}^\infty(X)}\Gamma(X,E) \cong
L_{{\cal C}^\infty(X)}\Big(\Gamma(X,E^*),\G(X)\Big)
$$
\et
In particular, this implies that generalized sections may be viewed as smooth sections with
generalized coefficients (in complete analogy to the distributional case). In addition,
for spaces of generalized tensor fields we have
\begin{itemize}
\item[]  $\G^r_s(X)\cong
L_{\G(X)}\Big(\G^0_1(X)^r,\G^1_0(X)^s;\G(X)\Big)$ as $\G(X)$-module. 
\item[] 
$\G^r_s(X)\cong L_{\CC^\infty(X)}\Big(\Om^1(X)^r,{\mathfrak X}(X)^s;\G(X)\Big)$
as $\CC^\infty(X)$-module. 
\end{itemize}
Contrary to the purely distributional picture where ill-defined products of distributions
have to be avoided carefully, our current setting allows unrestricted
application of multilinear operations like tensor product, wedge product, Lie derivatives
w.r.t.\ generalized vector fields, Poisson brackets, etc.

The relationship to the distributional setting is again governed by the  notion of
association: a generalized section $s\in\Ga_{\gs}(X,E)$ is called {\em associated to $w\in
\D'(X,E)$},
$s\approx w$, if for all $\mu\in\Ga_c(X,E^*\otimes\Vol(X))$ and one (hence every) 
representative $(s_\eps)_\eps$ of $s$
$$
        \lim\limits_{\eps\to 0}\int_X\,(s_\eps|\mu)\,=\,\langle w,\mu\rangle\,.
$$
Here, $(.|.)$ denotes the natural pairing
$$
\mbox{tr}_E\otimes \mbox{id}\mbox{: }
 (E\otimes E^*)\otimes \mbox{Vol}(X) \to (X\times \C)\otimes \mbox{Vol}(X)
= \mbox{Vol}(X)
$$
Stronger notions of association like $\approx_k$ are defined analogously to the scalar case.
Typically, multilinear operations on generalized sections display compatibility properties
with their distributional counterparts expressible in terms of association relations. E.g.,
if $\xi\in \gs^1_0(X)$ and $\xi\approx \eta\in {\D'}^1_0(X)$, $t\in \gs^r_s(X)$, $t\approx_\infty
u\in {\mathcal T}^r_s(X)$, then $L_\xi(t)\approx L_\eta(u)$.

Furthermore, classical theorems of smooth and distributional analysis (cf.\cite{marsden}) 
like the Poincar\'e lemma, Stokes' theorem, or the characterization of generalized vector fields
as derivations on generalized functions can be extended to the Colombeau setting (\cite{ndg, book}).
 
\section{Manifold-valued generalized functions} \label{mfval}
When applying generalized function techniques to problems of global analysis 
one inevitably encounters situations where a concept of generalized functions
defined on a manifold $X$ and taking values in another manifold is needed. 
Examples include flows of generalized vector
fields or geodesics of distributional spacetime metrics. Within
classical distribution theory, clearly no such concept is available. Colombeau
algebras on the other hand put more emphasis on the function-character
of the generalized functions (as opposed to the description as linear functionals
on spaces of test functions in the $\D'$-setting), which allows to develop an
appropriate theory in this framework. One main requirement with respect to such
a construction is that it be functorial. In particular, it must allow for 
unrestricted composition of generalized functions. In the local case, the problem
of composition of Colombeau functions was first addressed in \cite{AB}. The construction
suggested there formed the basis for the manifold case presented in \cite{gfvm,gfvm2}.
Since Colombeau functions by construction are localized on compact subsets of
their domain (in the sense that they are completely determined by the behavior of
their representatives on such sets, for small values of the regularization parameter),
in order to satisfy this requirement we have to single out representatives 
$(u_\eps)_\eps \in \cc^\infty(X,Y)^I$ which are {\em compactly bounded} (or {c-bounded}) 
in the following sense:
$$
\forall K\comp X\ \exists \eps_0>0\  \exists K'\comp Y \  \forall 
             \eps<\eps_0:\ u_\eps(K) \subseteq K'\,.
$$
Moderateness of nets $(u_\eps)_\eps \in \cc^\infty(X,Y)^I$, on the other hand, is formulated
using local charts. We thus arrive at the following definition:
\bd \label{modmapmf} 
The space $\esm[X,Y]$  
of compactly bounded 
(c-bounded)
moderate maps from $X$ to $Y$ is defined as the set
of all $(u_\eps)_\eps \in \cinfty(X,Y)^I$ such that
\begin{itemize}
  \item[(i)] $(u_\eps)_\eps$ is c-bounded
  \item[(ii)] 
   $\forall k\in\N$,
   for each chart 
   $(V,\vphi)$
   in $X$, each chart $(W,\psi)$ in $Y$, each $L\comp V$
   and each $L'\comp W$
   there exists $N\in \N$  with
   $$\sup\limits_{p\in L\cap u_\eps^{-1}(L')} \|D^{(k)}
   (\psi\circ u_\eps \circ \vphi^{-1})(\vphi(p))\| =O(\eps^{-N}).
   $$
\end{itemize}
\et
Note that the ``safety compact sets'' $L$ and $L'$ in this definition are needed
in order to control the potentially arbitrarily fast growth of chart diffeomorphisms
towards the boundary of their domains.

In the absence of a linear structure on the target space $Y$, we have to introduce
an equivalence relation in $\esm[X,Y]$ which precisely reduces to negligibility
of differences of representatives in the case $Y=\R^m$. We do this in a two step 
process. First, we assure that the distance between representatives as measured in
any Riemannian metric on $Y$ goes to zero. Growth conditions on
derivatives are then formulated in local charts:
\bd\label{equrel} Two elements $(u_\eps)_\eps$, $(v_\eps)_\eps$ of $\esm[X,Y]$ are 
called equivalent, $(u_\eps)_\eps \sim (v_\eps)_\eps$, 
if the following conditions are satisfied: 
\begin{itemize}
\item[(i)]  For all $K\comp X$, $\sup_{p\in K}d_h(u_\eps(p),v_\eps(p)) \to 0$ 
($\eps\to 0$)
for some (hence every) Riemannian metric $h$ on $Y$.
\item[(ii)] $\forall k\in \N_0\ \forall m\in \N$,
for each chart 
   $(V,\vphi)$
   in $X$, each chart $(W,\psi)$ in $Y$, each $L\comp V$
   and each $L'\comp W$:
$$
\sup\limits_{p\in L\cap u_\eps^{-1}(L')\cap v_\eps^{-1}(L')}\!\!\!\!\!\!\!\!\!\!\!\!\!\!\!\!\!\!\!
\|D^{(k)}(\psi\circ u_\eps\circ \vphi^{-1}
- \psi\circ v_\eps\circ \vphi^{-1})(\vphi(p))\|
=O(\eps^m). 
$$
\end{itemize}
\et
Finally, we define the space of Colombeau generalized functions defined on $X$ and taking values
in $Y$ as $\gs[X,Y]:=\esm[X,Y]/\sim$. Elements of $\gs[X,Y]$ typically model jump discontinuities,
whereas delta-type singularities are excluded by the c-boundedness of representatives (on the other
hand, it seems unclear anyways what a delta-type singularity should be in a manifold without additional
structure).

In analogy to (\ref{mgth}) one would expect that condition (ii) in \ref{equrel} need only
hold for $k=0$ in case $(u_\eps)_\eps$ is assumed to be moderate. It turns out, however, that a proof
of this fact cannot be carried along the lines of the local result (based in turn on a classical
argument by Landau, \cite{Lan}). Similarly, one would hope for a point value characterization
of elements of $\gs[X,Y]$. However, in the absence of an analogue to (\ref{mgth}) this seems
difficult to obtain.

The remedy for both problems lies in a nonlocal characterization of c-boundedness, moderateness
and equivalence (\cite{gfvm2}, Sec.\ 3). The key idea is to replace composition with charts in the
target space by composition with globally defined smooth functions.
\bp\label{propneu}
Let $(u_\eps)_\eps\in\CC^\infty(X,Y)^I$. The following conditions are equivalent
\begin{itemize}
\item[(i)] $(u_\eps)_\eps$ is c-bounded.
\item[(ii)] $(f\circ u_\eps)_\eps$ is c-bounded for all $f\in \CC^\infty(Y)$.
\item[(iii)] $(f\circ u_\eps)_\eps$ is moderate of order zero for all $f\in \CC^\infty(Y)$, i.e.,
\[
 \forall K\subset\subset X\ \exists N\in\N: \sup_{p\in K}|f\circ u_\eps(p)|=O(\eps^{-N})
\]
for all $f\in \CC^\infty(Y)$.
\item[(iv)] $(u_\eps(x_\eps))_\eps\in Y_c$ for all $(x_\eps)_\eps\in X_c$.
\end{itemize}
\et
Based on this result, moderateness can be characterized as follows:
\bp \label{emchar}
Let $(u_\eps)_\eps \in \cinfty(X,Y)^I$. Then $(u_\eps)_\eps \in \esm[X,Y]$ if and only if
$(f\circ u_\eps)_\eps \in \esm(X)$ for all $f\in \cinfty(Y)$. 
\et
Finally, concerning the equivalence relation $\sim$ on $\esm[X,Y]$ we obtain:
\bt \label{equivchar} 
Let $(u_\eps)_\eps$, $(v_\eps)_\eps \in \esm[X,Y]$. The following statements are
equivalent:
\begin{itemize}
\item[(i)] $(u_\eps)_\eps \sim (v_\eps)_\eps$.
\item[(ii)] For every
Riemannian metric $h$ on $Y$, every $m\in \N$ and every $K\comp X$,  
$$
\sup_{p\in K}d_h(u_\eps(p),v_\eps(p)) = O(\eps^m) \qquad (\eps\to 0)\,.
$$ 
\item[(iii)] $(f\circ u_\eps - f\circ v_\eps)_\eps \in \ns(X)$ for all $f\in \cinfty(Y)$.
\end{itemize}
\et
Since by \cite{gfvm}, Th.\ 2.14, condition (ii) in \ref{equivchar} is equivalent with 
conditions \ref{equrel} (i) and (ii) with $k=0$, we obtain the desired characterization
of $\sim$. This in turn provides the key building block in the proof of the following 
point value description of manifold-valued generalized functions:
\bt Let $u=[(u_\eps)_\eps]\in\gs[X,Y]$ and $\tilde p = [(p_\eps)_\eps]\in \tilde X_c$. Then
$u(\tilde p):=[(u_\eps(p_\eps))_\eps]$ is a well-defined element of $\tilde Y_c$. Moreover, 
$u,\, v \in \gs[X,Y]$ are equal if and only if their point values in each generalized point
agree. 
\et
Once this point value characterization is established, also the problem of composition
of generalized functions can be resolved (\cite{gfvm}, Th.\ 2.16 and \cite{gfvm2}, Th.\ 3.6):
\bt \label{generalcomp} Let 
$u=[(u_\eps)_\eps] \in \gs[X,Y]$, $v=[(v_\eps)_\eps] 
\in \gs[Y,Z]$. Then $v\circ u := [(v_\eps\circ u_\eps)_\eps]$ is a well-defined
element of $\gs[X,Z]$. 
\et
Although by the c-boundedness of representatives the ``worst'' singularities that can
be modelled  by elements of $\gs[X,Y]$ are jump discontinuities it is to be expected 
that derivatives (i.e., tangent
maps) of such generalized maps will behave $\delta$-like. We must therefore provide for
a concept of generalized vector bundle homomorphisms (containing such tangent maps as
special cases) with substantially less restrictive growth conditions in the vector 
components.

\bd \label{homgdef} For $E\to X$, $F\to Y$ vector bundles, $\esm^{\mathrm{VB}}[E,F]$ 
is the set of all $(u_\eps)_\eps$ $\in$ $\mathrm{Hom}(E,F)^I$ satisfying
\begin{itemize}
\item[(i)] $(\underline{u_\eps})_\eps \in \esm[X,Y]$.
\item[(ii)] $\forall k\in \N_0\
\forall (V,\Phi)$
vector bundle chart in $E$,
$\forall (W,\Psi)$ vector bundle chart in $F$,
$\forall L\comp V\
\forall L'\comp W\ \exists N\in \N\ \exists \eps_1>0\
\exists C>0$ with
$$
\|D^{(k)}
(u_{\eps \mathrm{\Psi}\mathrm{\Phi}}^{(2)}(\vphi(p)))\|
\le C\eps^{-N}
$$
for all $\eps<\eps_1$ and all $p\in L\cap\underline{u_\eps}^{-1}(L')$, with $\|\,.\,\|$ 
any matrix norm.
\end{itemize}
\et

\bd \label{homgequ}
$(u_\eps)_\eps$, $(v_\eps)_\eps \in \esm^{\mathrm{VB}}[E,F]$
are called $vb$-equivalent, $((u_\eps)_\eps \sim_{vb} (v_\eps)_\eps)$
if
\begin{itemize}
\item[(i)] $(\underline{u_\eps})_\eps \sim (\underline{v_\eps})_\eps$ in
$\esm[X,Y]$.
\item[(ii)] $\forall k\in \N_0\ \forall m\in \N\ \forall (V,\Phi)$
vector bundle chart in
$E$, $\forall (W,\Psi)$ vector bundle chart in $F$,
$\forall L\comp V\ \forall L'\comp W
\ \exists \eps_1>0\ \exists C>0$ such that:
$$
\|D^{(k)}(u_{\eps \mathrm{\Psi}\mathrm{\Phi}}^{(2)}
-v_{\eps \mathrm{\Psi}\mathrm{\Phi}}^{(2)})(\vphi(p))\|
\le C\eps^{m}
$$
for all $\eps<\eps_1$ and all $p\in L\cap\underline{u_\eps}^{-1}(L')
\cap\underline{v_\eps}^{-1}(L')$.
\end{itemize}
\et
We now set $\mathrm{Hom}_{\gs}[E,F] := \esm^{\mathrm{VB}}[E,F]\big/\sim_{vb}$.
For $u\in \mathrm{Hom}_{\gs}[E,F]$, $\underline{u} :=[(\underline{u}_\eps)_\eps]$
is a well-defined element of $\gs[X,Y]$ uniquely characterized by $\underline{u}
\circ\pi_X = \pi_Y\circ u$. The tangent map $Tu:=[(Tu_\eps)_\eps]$ of any 
$u\in\gs[X,Y]$ is then a well-defined element of $\mathrm{Hom}_{\gs}[TX,TY]$. 

Also in the context of generalized vector bundle homomorphisms a global characterization
of moderateness is available:
\bp \label{emvbchar} 
Let $(u_\eps)_\eps \in \mathrm{Hom}(E,F)^I$. Then
$(u_\eps)_\eps \in \esm^{VB}[E,F]$ if and only if
$(\hat f \circ u_\eps)_\eps \in \esm^{VB}(E,\R\times\R^{m'})$ 
              for all $\hat f \in \mathrm{Hom}(F,\R\times \R^{m'})$.
\et
and similar for $\sim_{vb}$ (\cite{gfvm2}, Prop.\ 4.1 and Th.\ 4.2). Based on these
results, appropriate point value descriptions of elements of $\mathrm{Hom}_{\gs}[TX,TY]$
can be derived. As a final ingredient, in Theorem \ref{th:mfflow} below we shall make use 
of the {\em hybrid} space $\gs^h[X,F]$ whose elements are defined on $X$ and take
values in $F$, c-bounded in the base component and moderate in the vector component (\cite{gprg,gfvm2}).
All of the above constructions are functorial (with compositions defined unrestrictedly).
We do not go into the details here (cf.\ \cite{gprg, gfvm2}) but instead
turn to another concept which is of relevance in applications to non-smooth 
pseudo-Riemannian geometry (cf.\ Sec.\ \ref{riemann}). Denote by
$$
\mathrm{Hom}_u(E,F):=\{v\in \mathrm{Hom}(E,F) | \underline{v} = u\}
$$
the space of generalized vector bundle homomorphisms over the generalized map $u$.
While in the smooth setting the corresponding space can trivially be endowed with
a vector space structure, 
the main obstruction in extending this property to the present
context is that, a priori, representatives 
$(v_\eps)_\eps$, $(v_\eps')_\eps$ of elements $v$, $v'$
of $\mathrm{Hom}_u(E,F)$ need not project onto the same representative $(u_\eps)_\eps$
of $u=\underline{v}=\underline{v'}\in\gs[X,Y]$, so that
simple fiberwise addition is in general not possible. The following result (\cite{gfvm2}, 
Prop.\ 5.7 and Cor.\ 5.8) remedies this problem:
\bp \label{fibertransplanting} 
Let $u=[(u_\eps)_\eps]\in \gs[X,Y]$ and $v\in \mathrm{Hom}_u(E,F)$. Then there
exists a representative $(v_\eps)_\eps$ of $v$ such that $\underline{v_\eps} = u_\eps$
for all $\eps\in I$. Consequently, $\mathrm{Hom}_u(E,F)$ is a vector space.
\et
To conclude this section let us have a look at the problem of determining the
flow of a generalized vector field $\xi\in \gs^1_0(X)$. We first note that in
the distributional setting already the notion of the flow of a distributional
vector field $\zeta$ is problematic, as it would have to denote a ``manifold-valued 
distribution''.
In \cite{marsden}, a regularization approach is used to cope
with this problem, by  introducing a c-bounded sequence of smooth vector fields
$\xi_\eps$ approximating $\zeta$. Each $\xi_\eps$ has a classical flow
$\Phi^\eps$ and under certain assumptions the assignment 
$\Psi=\lim_{\eps\to 0}\Phi^\eps$ allows to associate
a measurable flow $\Psi$ to the distributional vector field
$\zeta$. This approach is naturally related to the Colombeau picture, where
any $\xi = (\xi_\eps)_\eps \approx \zeta$ can be viewed as a regularization
of the distributional vector field $\zeta$. We first give a basic
existence and uniqueness result for flows of generalized vector fields
(\cite{flows}, Th.\ 3.6):
\bt\label{th:mfflow}
  Let $(X,h)$ be a complete Riemannian manifold and suppose that $\xi\in\gs^1_0(X)$ satisfies
  \begin{itemize}
  \item[(i)] $\xi=[(\xi_\eps)]$ with each $\xi_\eps$ globally bounded with respect to $h$.
  \item[(ii)] For each differential operator $P\in{\cal P}(X,TX)$ of first order 
  and each $K\comp X$: $\sup_{p\in K} \|(P\xi_\eps)|_p\|_h \le C|\log\eps|$
  (with $h$ any Riemannian metric).
  \end{itemize}
  Then there exists a unique generalized function $\Phi\in\gs[\R\times X,X]$,
  the {\em generalized flow of $\xi$}, such that
  \beas\label{eq:mfflow1}
   \frac{d}{dt}\Phi(t,x)&=&\xi(\Phi(t,x))\quad\mbox{ in }\G^h[\R\times X,TX] \\
   \Phi(0,.)&=&\mathrm{id}_{X} \quad\mbox{ in }\G[X,X]   \label{eq:mfflow2}\\
   \Phi(t+s,.)&=&\Phi(t,\Phi(s,.)) \quad\mbox{ in }\G[\R^2\times X,X]\,. 
                                                         \label{eq:mfflow3}
  \eeas 
\et
\bex 
 Let $X=T^2 = S^1\times S^1$ and $\xi=[(\xi_\eps)_\eps] \in \gs^1_0(X)$ with
 \[
  \xi_\eps(e^{i\alpha},e^{i\beta})=(e^{i\alpha},e^{i\beta};1,1-\rho_{\sigma(\eps)}(\alpha)).
 \]
 Here, $\rho$ is a test function with unit integral and $\sigma(\eps) = |\log(\eps)|^{-1}$. 
 Then since $X$ is compact, each $\xi_\eps$ possesses a
 global flow $\Phi^\eps$ and $\Phi:=[(\Phi^\eps)_\eps]$ $\in$ $\gs[\R\times X,X]$ is
 the unique generalized flow of $\xi$. $\Phi$ possesses a discontinuous pointwise limit $\Psi$,
 namely
 $$
  \Phi^\eps(t;e^{i\al},e^{i\beta})
  =\left(\begin{array}{l}
   e^{i(\alpha +t)}\\e^{i(\beta+t-\int\limits_\alpha^{\alpha+t}
   \rho_{\sigma(\eps)}(\gamma)\,d\gamma)}
  \end{array}\right)
  \,\to\,
  \left(\begin{array}{l}
   e^{i(\al+t)}\\e^{i(\beta+t-H(\alpha+t)+H(\alpha))}
  \end{array}\right),
 $$
 which satisfies the flow property  $\Psi_{s+t}=\Psi_s\circ\Psi_t$ for all $s,t\in\R$. 
\et
In general the question whether the unique generalized flow of a generalized vector field
possesses a limiting (measurable) flow is quite involved, cf.\ \cite{marsden, flows}.

\section{Generalized connections and non-smooth Riemannian geometry} \label{riemann}
Applications in general relativity have constituted one
of the main driving forces behind the development of non-smooth differential
geometry in the setting of Colombeau generalized functions (see \cite{vickersESI}). 
As an example, we consider so called {\em impulsive pp-waves}
(i.e., impulsive gravitational waves with parallel rays, cf.\ \cite{herbertgeo, geo}).
These are described by a distributional pseudo-Riemannian metric with line-element
\begin{equation}
  \label{ppmetric}
   ds^2=f(x,y)\delta(u)du^2-dudv+dx^2+dy^2\,.
\end{equation}
To extract physically relevant information from this spacetime metric one has to
be able to calculate curvature quantities and find solutions of the 
corresponding geodesic equations (determining the trajectories of particles
in the spacetime at hand). However, all of these operations are undefined within 
linear distribution theory: the former due to the nonlinear operations involved
in their calculation, the latter due to the lack of a concept of manifold-valued
distributions. On the other hand, as we have seen in the previous sections, 
algebras of generalized functions make available all the necessary tools to 
address these issues.
 
The following result forms the basis for the description of singular pseudo-Riemannian
metrics in the Colombeau framework (\cite{gprg}, Th.\ 3.1):
\bt \label{mainmetric}
For any generalized $(0,2)$-tensor $g\in \gs^0_2(X)$, the following are equivalent:
\begin{itemize}
\item[(i)] For each chart $(V_\al,\psi_\al)$ and each $\tilde p \in 
(\psi_\al(V_\al))^{\sim}_c$ the map
$g_\al(\tilde p): \gK^n\times\gK^n\to\gK^n$ is symmetric and nondegenerate.
\item[(ii)] $g:
\gs^1_0(X) \times \gs^1_0(X) \to \G(X)$ is symmetric and
$\mbox{det}(g)$ is invertible in $\G(X)$.
\item[(iii)] $\mbox{det}(g)$ is invertible in $\G(X)$ and for each relatively
compact open set $V\subseteq X$ there exists a representative $(g_\eps)_\eps$ of
$g$ and an $\eps_0 > 0$ such that $g_\eps|_V$ is a smooth pseudo-Riemannian
metric for all $\eps<\eps_0$.
\end{itemize}
\et
\bd\label{gmetric}
Let $g\in{\gs}^0_2(X)$ satisfy the conditions in \ref{mainmetric}. 
If, in addition, there exists $j\in \N_0$ such that 
the index of the $g_\eps$ as in \ref{mainmetric} (iii) equals $j$, we call
$g$ a generalized pseudo-Riemannian metric of index $j$ and $(X,g)$ a 
generalized pseudo-Riemannian manifold. If $j=1$ or $j=n-1$, $(X,g)$ is called a
generalized spacetime.
\et
It follows from finite-dimensional perturbation theory that the index so defined
does not depend on the chosen representative $(g_\eps)_\eps$ of $g$. 
With respect to applications, the most important characterization in Th.\
\ref{mainmetric} is (iii), as it guarantees that locally any generalized
metric has a representative consisting entirely of smooth pseudo-Riemannian 
metrics.

We note
first that the above way of modelling singular metrics is considerably more
flexible than the purely distributional approach: 
In \cite{marsden}, a distributional $(0,2)$-tensor field
$g\in{\D'}^0_2(X)$ is called nondegenerate if $g(\xi,\eta) = 0$ for all
$\eta\in \mathfrak{X}(X)$ implies $\xi=0\in\mathfrak{X}(X)$, while
in \cite{parker}, $g$ is called nondegenerate if it is nondegenerate 
(in the classical sense) off its singular support. The drawback of the
first definition is its ``nonlocality'', which is too weak to reproduce
the classical notion: e.g., $ds^2 = x^2\,dx^2$ is nondegenerate in this sense
although it is clearly singular at $x=0$. The second notion, on the other
hand, does not provide any restrictions on $g$ at its points of singularity.

Since $\gs(X)$ is an algebra, all curvature quantities (Riemann-tensor, Ricci-
and Einstein tensor, \dots) of a generalized metric can be calculated unrestrictedly.
Moreover, in parallel to the smooth setting, we may develop a generalized 
pseudo-Riemannian geometry based on the above notions. Our first basic result towards that
goal is the following (\cite{gprg}, Prop.\ 3.9):
\bp \label{mlemma}
Let $(X,g)$ be a generalized pseudo-Riemannian manifold.
\begin{itemize}
\item[(i)] $g$ is non-degenerate in the following sense:
if $\xi\in \gs^1_0(X)$ and $g(\xi,\eta)=0$ $\forall \eta\in{\gs}^1_0(X)$, then $\xi=0$.
\item[(ii)] $ g$ induces a $\gs(X)$-linear isomorphism 
${\gs}^1_0(X)\to{\gs}^0_1(X)$
by $\xi\mapsto  g(\xi,\,.\,)$.
\end{itemize}
\et
The isomorphism in (ii) can naturally be extended to higher order tensor fields, so that,
as in the smooth case, generalized metrics can be used to raise and lower indices.
\bd\label{gcon}\ 
A {\em generalized connection $\hat D$} on $X$ is a map
${\gs}^1_0(X)\times{\gs}^1_0(X)\to{\gs}^1_0(X)$ satisfying 
\begin{itemize}
\item[(D1)] $\hat D_\xi \eta$ is $\gR$-linear in $\eta$.
\item[(D2)] $\hat D_\xi \eta$ is $\gs(X)$-linear in $\xi$.
\item[(D3)] $\hat D_\xi(u\eta)=u\,\hat D_\xi \eta+\xi(u)\eta$ for all $u\in\gs(X)$.
\end{itemize}
\et
With this notion we have the following
{\em Fundamental Lemma of pseudo-Riemannian geometry} (\cite{gprg}, Th. 5.2):
\bt 
On each generalized pseudo-Riemannian manifold $(X, g)$ there exists a unique
generalized Levi-Civita connection $\hat D$ 
such that for all $\xi,\,\eta,\,\zeta$ in ${\gs}^1_0(X)$:
\begin{itemize}
\item[(D4)] $[\xi,\eta]=\hat D_\xi \eta-\hat D_\eta\xi$ and
\item[(D5)] $\xi\, g(\eta,\zeta)= g(\hat D_\xi \eta,\zeta)+  g(\eta,\hat D_\xi \zeta)$
\end{itemize}
\et
Suppose now that $\ga\in \gs[J,X]$ is a generalized curve in $X$ defined on some interval
$J\subseteq \R$. Using a representative $(g_\eps)_\eps$ as in Th.\ \ref{mainmetric} (iii)
we may componentwise define an induced covariant derivative $\xi \mapsto \xi'$ on the space 
$\mathfrak{X}_{\gs}(u) := \{\xi\in\gsh[X,TY] \mid 
\underline{\xi} = u \}$ of generalized vector fields on $\ga$. Its basic properties are
summarized in the following result (\cite{gprg}, Prop.\ 5.6 and \cite{gfvm2}, Sec.\ 5):
\bt
Let $(X,g)$ be a generalized pseudo-Riemannian manifold and let $\gamma \in \gs[J,X]$. Then 
\begin{itemize}
\item[(i)] $(\tilde r \xi_1 + \tilde s \xi_2)' = \tilde r \xi_1' + \tilde s \xi_2'$ 
$\qquad (\tilde r,\, \tilde s \in \ks, \, \xi_1,\, \xi_2 \in \mathfrak{X}_{\gs}(\gamma))$.
\item[(ii)] $(u\xi)' = \frac{du}{dt}\xi + u\xi'$ $\qquad (u\in \gs(J),\,
\xi \in \mathfrak{X}_{\gs}(\gamma))$.
\item[(iii)] $(\xi\circ \gamma)' = \hat D_{\gamma'(.)}\xi$ 
\hspace{1em} in $\mathfrak{X}_{\gs}(\gamma)$ $(\xi \in {\gs}^1_0(X))$.
\item[(iv)] $\frac{d}{dt}g(\xi,\eta) = g(\xi',\eta) + g(\xi,\eta')$
\qquad $(\xi,\, \eta \in \mathfrak{X}_{\gs}(\gamma))$.
\end{itemize}
\et
Note in particular that property (iv) only makes sense due to Proposition \ref{fibertransplanting}.
Now that we have induced covariant derivatives at our disposal we may as  in the smooth case (and 
contrary to the distributional setting) give the following definition:
\bd A curve  $\ga\in\gs[J,X]$ in a generalized pseudo-Riemannian manifold is called
geodesic if $\ga''=0$. Here $\ga''$ is the induced covariant derivative of the
velocity vector field $\ga'$ of $\ga$.
\et 
Locally, therefore, the determination of the geodesics of a given singular metric
amounts to the solution of a system of ordinary differential equations in the
Colombeau setting. This program has been carried out for our first example 
(\ref{ppmetric}) in \cite{geo,geo2}. Using a generic regularization procedure for
the delta-term in (\ref{ppmetric}), the resulting system is uniquely solvable
in $\gs[\R,X]$. Moreover, for $\eps\to 0$ (i.e., in the sense of association) 
this unique solution displays the physically expected behavior of broken, 
refracted straight lines as geodesics. 

As a further aspect of the spacetime (\ref{ppmetric}) we note that its analysis naturally
leads to the concept of manifold-valued generalized functions: In \cite{penrose2},
R.\ Penrose introduced a discontinuous coordinate transformation $T$ that formally
transforms the distributional metric (\ref{ppmetric}) into a continuous form.
Although the two forms of the metric are physically equivalent (in the sense that
they have the same geodesics), the transformation relating them is clearly ill-defined 
in the distributional picture. In \cite{penrose}, however, $T$ was identified as an
element $[(T_\eps)_\eps]$ of $\gs[X,X]$ with each $T_\eps$ a diffeomorphism.
In this sense $T$ itself may be considered a ``discontinuous diffeomorphism''.

Recently, generalized pseudo-Riemannian geometry in the sense of the present section
has been identified as a special case of an encompassing theory of generalized connections
on fiber bundles. For this theory as well as for first applications to singular solutions
of Yang-Mills equations we refer to \cite{conn}.


\begin{thebibliography}{10}

\bibitem{AB}
{Aragona, J., Biagioni, H.~A.}
\newblock Intrinsic definition of the {C}olombeau algebra of generalized
  functions.
\newblock {\em Analysis Mathematica}, {\bf 17}:75--132, 1991.

\bibitem{herbertgeo}
{Balasin, H.}
\newblock Geodesics for impulsive gravitational waves and the multiplication of
  distributions.
\newblock {\em Class.~Quant.~Grav.}, {\bf 14}:455--462, 1997.

\bibitem{CPort}
{Colombeau, J.~F.}
\newblock New generalized functions. {M}ultiplication of distributions.
  {P}hysical applications. {C}ontribution of {J}. {S}ebastiao e {S}ilva.
\newblock {\em Port. Math.}, {\bf 41}:57--69, 1982.

\bibitem{CJMAA}
{Colombeau, J.~F.}
\newblock A multiplication of distributions.
\newblock {\em J. Math. Anal. Appl.}, {\bf 94}:96--115, 1983.

\bibitem{CCR}
{Colombeau, J.~F.}
\newblock Une multiplication g\'en\'erale des distributions.
\newblock {\em C. R. Acad. Sci. Paris, S\'er. I}, {\bf 296}:357--360, 1983.

\bibitem{c1}
{Colombeau, J.~F.}
\newblock {\em New Generalized Functions and Multiplication of Distributions}.
\newblock North Holland, Amsterdam, 1984.

\bibitem{c2}
{Colombeau, J.~F.}
\newblock {\em Elementary Introduction to New Generalized Functions}.
\newblock North Holland, Amsterdam, 1985.

\bibitem{DP}
{Dapi\'c, N., Pilipovi\'c, S.}
\newblock Microlocal analysis of {C}olombeau's generalized functions on a
  manifold.
\newblock {\em Indag.~Math.~(N.S.)}, {\bf 7}:293--309, 1996.

\bibitem{deR}
{De Rham, G.}
\newblock {\em Differentiable Manifolds}, volume~{\bf 266} of {\em Grundlehren
  der mathe\-mati\-schen Wissenschaften}.
\newblock Springer, Berlin, 1984.

\bibitem{RD}
{De Roever, J.~W., Damsma, M.}
\newblock Colombeau algebras on a {${\cal C}^\infty$}-manifold.
\newblock {\em {Indag.~Mathem., N.S.}}, {\bf 2}(3), 1991.

\bibitem{gt}
{Geroch, R., Traschen, J.}
\newblock Strings and other distributional sources in general relativity.
\newblock {\em Phys.~Rev.~D}, {\bf 36}(4):1017--1031, 1987.

\bibitem{found}
{Grosser, M., Farkas, E., Kunzinger, M., Steinbauer, R.}
\newblock On the foundations of nonlinear generalized functions {I}, {II}.
\newblock {\em Mem. Amer. Math. Soc.}, {\bf 153}(729), 2001.

\bibitem{book}
{Grosser, M., Kunzinger, M., Oberguggenberger, M., Steinbauer, R.}
\newblock {\em Geometric Theory of Generalized Functions}, volume 537 of {\em
  Mathematics and its Applications 537}.
\newblock Kluwer Academic Publishers, Dordrecht, 2001.

\bibitem{vim}
{Grosser, M., Kunzinger, M., Steinbauer, R., Vickers, J.}
\newblock A global theory of algebras of generalized functions.
\newblock {\em Adv. Math.}, 166:179--206, 2002.

\bibitem{hermannbook}
{Hermann, R.}
\newblock {\em {C-O-R} Generalized Functions, Current Algebras, and Control},
  volume~{\bf 30} of {\em Interdisciplinary Mathematics}.
\newblock Math Sci Press, 1994.

\bibitem{gfvm}
{Kunzinger, M.}
\newblock Generalized functions valued in a smooth manifold.
\newblock {\em Monatsh.\ Math.}, {\bf 137}:31--49, 2002.

\bibitem{flows}
{Kunzinger, M., Oberguggenberger, M., Steinbauer, R., Vickers, J.}
\newblock Generalized flows and singular odes on differentiable manifolds.
\newblock {\em Acta Appl. Math.}, to appear, 2004.

\bibitem{penrose}
{Kunzinger, M., Steinbauer, R.}
\newblock A note on the {P}enrose junction conditions.
\newblock {\em Class. Quant. Grav.}, {\bf 16}:1255--1264, 1999.

\bibitem{geo2}
{Kunzinger, M., Steinbauer, R.}
\newblock A rigorous solution concept for geodesic and geodesic deviation
  equations in impulsive gravitational waves.
\newblock {\em J. Math. Phys.}, {\bf 40}(3):1479--1489, 1999.

\bibitem{ndg}
{Kunzinger, M., Steinbauer, R.}
\newblock Foundations of a nonlinear distributional geometry.
\newblock {\em Acta Appl. Math.}, 71:179--206, 2002.

\bibitem{gprg}
{Kunzinger, M., Steinbauer, R.}
\newblock Generalized {pseudo-}{R}iemannian geometry.
\newblock {\em Trans. Amer. Math. Soc.}, 354(10):4179--4199, 2002.


\bibitem{gfvm2}
{Kunzinger, M., Steinbauer, R., Vickers, J.}
\newblock Intrinsic characterization of manifold-valued generalized functions.
\newblock {\em Proc.\ London Math.\ Soc.}, {\bf 87}(2):451--470, 2003.

\bibitem{conn}
{Kunzinger, M., Steinbauer, R., Vickers, J.}
\newblock Generalized connections and curvature.
\newblock {\em Preprint}, 2003.


\bibitem{Lan}
{Landau, E.}
\newblock Einige {U}ngleichungen f\"ur zweimal differentiierbare {F}unktionen.
\newblock {\em Proc.~London Math.~Soc. Ser.~2}, {\bf 13}:43--49, 1913--1914.

\bibitem{marsden}
{Marsden, J.~E.}
\newblock Generalized {H}amiltonian mechanics.
\newblock {\em Arch.~Rat.~Mech.~Anal.}, {\bf 28}(4):323--361, 1968.

\bibitem{marsden2}
{Marsden, J. E.}
\newblock Non-smooth geodesic flows and classical mechanics.
\newblock {\em Canad. Math. Bull.}, 12:209--212, 1969.

\bibitem{MObook}
{Oberguggenberger, M.}
\newblock {\em Multiplication of Distributions and Applications to Partial
  Differential Equations}, volume~{\bf 259} of {\em Pitman Research Notes in
  Mathematics}.
\newblock Longman, Harlow, U.K., 1992.

\bibitem{point}
{Oberguggenberger, M., Kunzinger, M.}
\newblock Characterization of {C}olombeau generalized functions by their
  pointvalues.
\newblock {\em Math. Nachr.}, {\bf 203}:147--157, 1999.

\bibitem{OPS}
{Oberguggenberger, M., Pilipovic, S., Scarpalezos, D.}
\newblock Local properties of {C}olombeau generalized functions.
\newblock {\em Math. Nachr.}, {\bf 256}:88--99, 2003.

\bibitem{parker}
{Parker, P.}
\newblock Distributional geometry.
\newblock {\em J.~Math.~Phys.}, {\bf 20}(7):1423--1426, 1979.

\bibitem{penrose2}
{Penrose, R.}
\newblock The geometry of impulsive gravitational waves.
\newblock In L.~O'Raifeartaigh, editor, {\em General Relativity, Papers in
  Honour of {J}.\,{L}. {S}ynge}, pages 101--115. Clarendon Press, Oxford, 1972.

\bibitem{Schw}
{Schwartz, L.}
\newblock Sur l'impossibilit\'e de la multiplication des distributions.
\newblock {\em C.~R.~Acad.~Sci.~Paris}, {\bf 239}:847--848, 1954.

\bibitem{geo}
{Steinbauer, R.}
\newblock Geodesics and geodesic deviation for impulsive gravitational waves.
\newblock {\em J.~Math.~Phys.}, {\bf 39}(4):2201--2212, 1998.

\bibitem{vickersESI}
{Vickers, J.~A.}
\newblock Nonlinear generalized functions in general relativity.
\newblock In {Grosser, M., Hörmann, G., Kunzinger, M., Oberguggenberger, M.},
  Eds., {\em Nonlinear Theory of Generalized Functions}, 275--290, CRC, Boca
  Raton 1999.

\bibitem{wilson}
{Wilson, J.~P.}
\newblock Distributional curvature of time dependent cosmic strings.
\newblock {\em Class.~Quantum Grav.}, {\bf 14}:3337--3351, 1997.

\end{thebibliography}
\end{document}